\documentclass{amsart}
\usepackage{epsf}
\def\hepsffile{\leavevmode\epsffile}

\theoremstyle{plain}
\newtheorem{thm}{Theorem}[subsection]

\newtheorem{lem}[thm]{Lemma}

\newtheorem{prop}[thm]{Proposition}

\theoremstyle{definition}
\newtheorem{defin}[thm]{Definition}

\newtheorem{emf}[thm]{}
\newtheorem{rem}[thm]{Remark}

\def\aslk{\protect\operatorname{aslk}}

\def\Im{\protect\operatorname{Im}}

\def\Int{\protect\operatorname{Int}}

\def\slk{\protect\operatorname{slk}}

\def\spin{\protect\operatorname{spin}}

\def\slk{\protect\operatorname{slk}}

\def\pr{\protect\operatorname{pr}}

\def\wt{\widetilde}

\def\C{{\mathbb C}}
\def\Z{{\mathbb Z}}

\def\R{{\mathbb R}}
\def\N{{\mathbb N}}

\def\1{\hbox{\rm\rlap {1}\hskip.03in{\rom I}}}
\def\Bbbone{{\rm1\mathchoice{\kern-0.25em}{\kern-0.25em}
	{\kern-0.2em}{\kern-0.2em}I}}

\def\p{\partial}

\begin{document}
\hyphenation{Ca-m-po}
\title[Framed knots in $3$-manifolds and affine self-linking numbers]
{Framed knots in $3$-manifolds and affine self-linking numbers}
\author[V.~Chernov (Tchernov)]{Vladimir Chernov (Tchernov)}
\address{Department of Mathematics, 6188 Bradley Hall, Dartmouth College, Hanover, NH 03755, USA}
\email{Vladimir.Chernov@dartmouth.edu}
\begin{abstract}
The number $|K|$ of non-isotopic framed knots that correspond to a given unframed knot $K\subset S^3$ is infinite. This follows from the existence of 
the self-linking number $\slk$ of a zerohomologous framed knot. We use the approach of Vassiliev-Goussarov invariants to construct ``affine self-linking numbers'' that are extensions of $\slk$ to the case of  nonzerohomologous framed knots. 

As a corollary we get that $|K|=\infty$ for all  knots in an oriented (not necessarily compact) $3$-manifold $M$ that is not realizable as a connected sum $(S^1\times S^2) \# M'$. This result for compact manifolds was first stated by Hoste and Przytycki. They referred to the works of McCullough for the idea of the proof, 
however to the best of our knowledge the proof of this fundamental fact was not given in literature. Our proof is based on different ideas.
For $M=(S^1\times S^2) \# M'$  we construct $K$ in $M$ such that $|K|=2\neq \infty$.
\end{abstract}
 
\maketitle

\leftline {\em \Small 2000 Mathematics Subject Classification.
Primary: 57M27}

\leftline{\em \Small Keywords: framed knots, self-linking and linking numbers, finite order
invariants}

\section{Introduction}\label{Introduction} 
We work in the smooth category. Throughout this paper $M$ is a smooth 
oriented connected $3$-dimensional Riemannian manifold, unless the opposite is explicitly stated.

A {\em curve\/} in a manifold $M$ is an immersion of $S^1$ into
$M$. A {\em framed curve\/} in $M$ is a curve equipped with a continuous unit
normal vector field.

A {\em knot (resp. framed knot)\/} 
in $M$ 
is a smooth embedding (resp. framed smooth embedding) of $S^1$. An {\em isotopy} of ordinary (framed) knots is a path in the space of 
ordinary (framed) smoothly embedded curves.

For an unframed knot $K$ in $M$ we denote by $|K|\in \{\N \sqcup \infty\}$ the number of isotopy classes of framed knots that correspond to the isotopy class of $K$ when one forgets the framing.

For a framed knot $K_f$ and $i\in \Z$ we denote by $K_f^i$ the framed knot obtained by addition  of $i$ extra positive twists to the framing of $K_f$, if $i\geq 0$; and by addition of $|i|$ extra negative twists to the framing of $K_f$, if $i<0$. It is clear that every isotopy class of framed knots corresponding to $K$ is realizable as one of $K_f^i$. However what could happen (and actually does happen for knots in some manifolds) 
is that $K_f^i$ and $K_f^j$ are isotopic for $i\neq j$. In this case 
the isotopy classes of $K_f^{i+n}$ and of $K_f^{j+n}$ are the same for every $n\in \Z$, and hence $|K|$ is finite.

It is well-known that if $M=S^3$, then $|K|=\infty$, and intuitively one expects that $|K|=\infty$ for all $M$.
For zerohomologous $K$ this follows from the existence of the self-linking number $\slk (K_f)$ of a zero-homologous framed $K_f$ defined as the intersection of the infinitesimal shift of $K_f$ along its framing with an oriented surface $S$ bounded by $K_f$. It is easy to verify that self-linking of a zero-homologous framed knot is indeed well-defined and in particular does not depend on the choice of the oriented surface $S$.
It is also clear that this definition can not be made unless the knot is zero-homologous. (The work of U.~Kaiser~\cite{Kaiserbook} allows one to extend the homology definition of self-linking to the case where the knot is homologous into the boundary of the ambient manifold.) 

The self-liking number $\slk$ can also be defined as a Vassiliev-Goussarov invariant satisfying certain axioms. We show that the last definition (and certain modifications of it) make sense for a vast collection of nonzerohomologous framed knots. We called the resulting invariants {\em affine self-linking numbers\/} $\aslk$. 
As a corollary we get that $|K|=\infty$ for all (not necessarily zero-homologous) $K$ in a not necessarily compact $M$ such that $M$ is not realizable as a connected sum $M'\#(S^1\times S^2)$. This result (for compact manifolds) was first stated by Hoste and Przytycki~\cite{HostePrzytycki}. 
They referred to the work~\cite{McCullough} of McCullough on mapping 
class groups of 3-manifolds for the idea of the proof of 
this fact. However to the best of our knowledge the proof  
was not given in the literature. The proof we provide is 
based on the ideas and methods different from the ones 
Hoste and Przytycki had in mind. 

For $M=M'\#(S^1\times S^2)$ the examples of $K$ such that $|K|\neq \infty$ were first constructed by Hoste and Przytycki~\cite{HostePrzytycki}, and are also discussed in this paper.

It is well-known that $M$ is realizable as a connected sum $(S^1\times S^2)\#M'$ if and only if $M$ contains an embedded sphere that does not separate $M$ into two parts. This allows one to reformulate the results of the paper in terms of non-separating spheres. 

The famous Bennequin inequality provides restrictions on the self-linking number of a zero-homologous framed knot realizable by Legendrian and transverse knots in a tight contact manifold. We believe that the affine self-linking numbers constructed in this paper might provide the necessary ingredient for the generalization of the Bennequin inequality to nonzerohomologous knots. What one should probably hope to prove is that only framed knots with an affine self-linking number less than a certain constant can be realized by Legendrian and transverse knots in a tight contact manifold.

These results of the author most of which appeared as a preprint~\cite{Chernovframedpreprint} are philosophically similar to the later joint results of the author and Yu.~Rudyak~\cite{ChernovRudyakgenerallinking} showing that ``affine linking numbers'' that generalize the classical linking numbers can be defined for a vast collection of nonzerohomologous submanifolds $N_1$ and $N_2$ of the ambient manifold $M$ such that $\dim N_1+\dim N_2+1=\dim M$. Affine linking numbers have interesting applications to causality and mathematical physics, 
see~\cite{ChernovRudyakWinding} and~\cite{ChernovRudyakCausality}. 

The results of this paper do not follow from our results on affine linking between two submanifolds~\cite{ChernovRudyakgenerallinking}, and in many cases affine self-linking numbers are harder to construct than affine linking numbers.
As it is shown in our works with Rudyak~\cite{ChernovRudyakgenerallinking} and~\cite{ChernovRudyakGerstenhaber} the obstruction for the existence of the ``affine linking number" is the image of the 
generalized string homology Lie bracket of Chas and Sullivan~\cite{ChasSullivan} introduced by us in~\cite{ChernovRudyakGerstenhaber}. (This Lie bracket also generalizes~\cite{ChernovRudyakgenerallinking} the Goldman Lie bracket~\cite{Goldman} of free loops on surfaces.)
This paper illustrates that, similar to linking, self-linking can be defined for many nonzerohomologous framed submanifolds via the approach of Vassiliev-Goussarov invariants, and that self-linking has nice applications in topology. However the general theory of affine self-linking invariants of framed submanifolds still needs to be developed.

\section{Main results}

\begin{defin}[Vassiliev-Goussarov invariants]
A {\em singular (framed)\/} knot with $n$ double points is a curve (framed curve)
in $M$ whose only singularities are $n$ transverse double points.
An {\em isotopy\/} of a singular (framed) knot 
with $n$ double points is a path in the space of singular (framed) knots with
$n$ double points under which the preimages of the double points on $S^1$
change continuously.

For an Abelian group $\mathcal A$ an $\mathcal A$-valued (framed) knot invariant is an $\mathcal A$-valued function on the set of  isotopy classes of (framed) knots.

A transverse double point $t$ of a singular knot can be resolved in two 
essentially different ways. We say that a resolution of a double point is
positive (resp. negative) if the tangent vector to the
first strand, the tangent vector to the second strand, and the vector from
the second strand to the first form the positive $3$-frame. (This does 
not depend on the order of the strands).

A singular (framed) knot $K$ with $(n+1)$ 
transverse double points
admits $2^{n+1}$ possible resolutions of the double points. The sign of the resolution 
is put to be $+$ if the number of negatively resolved double points is even, and
it is put to be $-$ otherwise. 
Let $x$ be an $\mathcal A$-valued invariant of (framed) knots. The invariant $x$ is said to be of {\em finite
order (or Vassiliev, or Vassiliev-Goussarov invariant\/}) if there exists a nonnegative 
integer $n$ such that for any singular knot $K_s$ with $(n+1)$
transverse double points the sum (with the signs defined above) of the values of $x$ on the nonsingular
knots obtained by the $2^{n+1}$ resolutions of the double points is zero. 
An invariant is said to be of order not greater than $n$ (of order $\leq n$) if $n$
can be chosen as the integer in the definition above. The group of $\mathcal
A$-valued finite order invariants has an increasing filtration by the
subgroups of the invariants of order $\leq n$.
\end{defin}

It is easy to verify that the {\em self-linking invariant\/} of a zero-homologous framed knot is a $\Z$-valued order $\leq 1$ Vassiliev-Goussarov invariant that increases by two under every positive passage through a transverse double point. 
(The sign of a  passage through a double point is the sign of the resolution of the double point that happens during the passage.) 
In fact the above properties define the self-linking number up to the choice of its value on one framed knot in each connected component of the space of framed curves. 

Recall that a manifold $M$ is said to be
{\em irreducible\/} if every two-sphere embedded into $M$ bounds a ball. It is well-known that every closed oriented $3$-manifold that is not realizable as a connected sum $(S^1\times S^2)\#M'$ admits a decomposition into a connected sum of irreducible manifolds. 
In fact this decomposition is unique up to the permutation of the summands and additions of $S^3$. Recall also that a manifold is {\em prime} if it can not be decomposed as a nontrivial connected sum, and that $(S^1\times S^2)$ is the 
only oriented closed prime manifold that is not irreducible.

\begin{thm}\label{aslk} Let $M$ be a closed oriented $3$-manifold realizable as a connected sum of irreducible $3$-manifolds $\#_{j\in J} M_j$, and let $\mathcal F'$ be a connected component of the space of framed curves in $M$. 
\begin{description}
\item[1] If framed curves from $\mathcal F'$  are not homotopic 
to a curve  contained in one of the irreducible pieces $M_j$, then there exists a $\Z$-valued ``affine self-linking invariant'' $\aslk$ of framed knot from $\mathcal F'$ such that it increases by two under every positive passage through a transverse double point. 
\item[2] If curves from $\mathcal F'$ are homotopic to a curve contained in one of the irreducible pieces $M_j$, then there exists a $\Z$-valued ``affine self-linking invariant'' $\wt \aslk$ of framed knot from $\mathcal F'$ such that: 
\begin{description}
\item[a] it increases by two under every positive passage through a transverse double point, provided that one of the two loops of the singular knot
separated by the double point is contractible; 
\item[b] it does not change under other passages through a double point.
\end{description}
\end{description}
\end{thm}

For the Proof of Theorem~\ref{aslk} see Section~\ref{Proofaslk}.

\begin{rem} 
Invariants $\aslk$ and $\wt \aslk$ defined in the Theorem are clearly Vassiliev  invariants of order $\leq 1$. It is easy to see from the proof of Theorem~\ref{aslk} that they are defined uniquely up to an additive constant which is the value of an invariant
on a preferred knot from $\mathcal F'$. This ambiguity in the choice of the value on a preferred knot is similar to the ambiguity 
in the choice of the zero vector in the affine vector space, and it is the reason for the word affine in the name of the invariants.

As it was explained before, the definition of the $\aslk$ invariant coincides with the Vassiliev-Goussarov definition of the classical self-linking number $\slk$ of a zero-homologous framed knot. 
Since $\pi_1(S^3)=1$, it is easy to see that the invariant $\wt\aslk$ coincides (up to an additive constant) with the classical self-linking number of framed knots in $S^3$. Thus both invariants $\aslk$ and $\wt\aslk$ are natural generalizations of the self-linking number.

A closed irreducible manifold is {\em atoroidal\/} if it does not admit {\em essential mappings\/} $\mu:S^1\times S^1\to M$, i.e. mappings such that $\mu_*:\pi_1(S^1\times S^1)\to \pi_1(M)$ is injective. As it is essentially shown in the first version of~\cite{ChernovLegendrianpreprint}, Theorem 3.0.8, and in~\cite{Chernovframedpreprint}, the $\aslk$ invariant exists also for all knots in an irreducible closed atoroidal $3$-manifold $M$. (The same result can be easily concluded from the slightly later independent works of 
U.~Kaiser~\cite{Kaisermodule} and~\cite{Kaisermodulenotframed}.)
As we show in~\ref{nonessentialaslk} and~\ref{essentialorientable} the same in true if the irreducible component $M_j$ of $M$ into which $K$ is homotopic is atoroidal; or if $M_j$ is not one of the Seifert-fibered manifolds over $S^2$ corresponding to Euclidean triangle groups and the characteristic submanifold of $M_j$ does not contain components  that are Seifert-fibered over a nonorientable surface.
\end{rem}

\begin{emf}{\em Affine self-linking invariants for framed knots in compact $M$.\/} As it is shown in Lemma~\ref{Matveevlemma} told to us by S.~Matveev~\cite{Matveev}, every compact oriented  $M$ that does not contain nonseparating spheres can be included $i:M\to \overline M$ into a closed oriented  $\overline M$ that also does not contain embedded nonseparating spheres. (Recall that $M$ contains an embedded  nonseparating sphere if and only if $M=(S^1\times S^2)\#M'$, for some $M'$.)
This allows one to define affine self-linking numbers
for framed knots in compact $M$. 

However $i$ does not necessarily respect the decomposition into irreducible submanifolds, since $M$ could contain separating disks, and $i_*:\pi_1(M)\to \pi_1(\overline M)$ is not always injective.
This approach allows to define at least one of $\wt \aslk$ and $\aslk$ for a framed 
knot in an oriented compact $M$ that does not contain nonseparating spheres. It can occur however that the inclusion construction defines the $\aslk$ invariant (rather than $\wt \aslk$) for a knot homotopic into an irreducible summand of $M$ (but not of $\overline M$). Also the $\wt \aslk$ invariant of a knot in $M$ defined via $i$ might change under passages through double points with one loop of the singular knot contractible in $\overline M$, but not necessarily in $M$. 

Luckily in many cases it is possible to show the existence of $\aslk$ and $\wt \aslk$ directly. Namely in the proof of statement {\bf 1} of Theorem~\ref{aslk} it suffices to assume that $M$ is compact (and not closed), see~\ref{aslkcompactnotclosed}. 
In the proof of statement {\bf 2} of Theorem~\ref{aslk} it suffices to assume that
the prime summand $M_j$ of $M$ into which knots from $\mathcal F'$ are homotopic either does not admit essential mappings of $S^1\times S^1$ or that it is Haken, see~\ref{wtaslkcompacthaken}.
In particular $\wt \aslk$ exists if every essential mappings of a torus into $M_j$ is homotopic into  $\p M_j$, a conclusion that is essentially contained in the work of U.~Kaiser~\cite{Kaisermodule} for irreducible $M$. 
\end{emf}

The affine self-linking invariants of Theorem~\ref{aslk} are essential in the proof of the following Theorem.

\begin{thm}\label{framed} 
Let $M$ be a (not necessarily compact) oriented $3$-manifold such that it is not realizable as a connected sum $M'\#(S^1\times S^2)$, then  $|K|=\infty$ for all (not necessarily zero-homologous) unframed knots $K$ in $M$.
\end{thm}

For the proof of Theorem~\ref{framed} see Section~\ref{proofframed}.

This result (for compact manifolds) was first stated by Hoste and Przytycki~\cite{HostePrzytycki}. 
They referred to the work~\cite{McCullough} of McCullough on mapping 
class groups of 3-manifolds for the idea of the proof of 
this fact. However to the best of our knowledge the proof  
was not given in the literature. The proof we provide is 
based on the ideas and methods different from the ones 
Hoste and Przytycki had in mind.
(Partial cases of this result of the author prior to the preprint~\cite{Chernovframedpreprint} appeared in~\cite{ChernovLegendrianpreprint}, Theorems 3.0.6, 3.0.9, 
and Remark
3.0.10, where the Theorem 
was proved  for all knots in closed orientable 
hyperbolic $3$-manifolds and in the orientable total spaces
of locally-trivial $S^1$-bundles over a surface $F\neq S^2, \R P^2$.)

\begin{emf}{\em Framed knots in non-orientable $3$-manifolds.\/}
Straightforward geometric considerations show that if $M$ is a non-orientable
$3$-manifold, then $|K|=2$ for every knot that realizes an orientation
reversing loop in $M$.

It is also clear that if $K$ realizes an orientation preserving loop in $M$ and $|K|\neq \infty$, then $|\tilde K|\neq \infty$ for the lifting $\tilde K$ of $K$ to the total space $\tilde M$ of the orientation double cover $p:\tilde M\rightarrow M$. 
Thus if $M$ is a non-orientable (not necessarily compact) connected 
$3$-manifold such that $\tilde M$ is not realizable as $(S^1\times S^2)\#\tilde M'$, then $|K|=\infty$ for every $K$ realizing an orientation 
preserving loop.
\end{emf}

The following Theorem shows that if $M=M'\#(S^1\times S^2)$, then 
there exist $K\subset M$ such that $|K|$ is finite. Examples of this sort were previously described by Hoste and Przytycki~\cite{HostePrzytycki} and later in the works of the author, see~\cite{ChernovLegendrianpreprint} Theorem 3.1.2.a
and~\cite{ChernovLegendrian}.

\begin{thm}\label{Example}
Let $M$ be an oriented (not necessarily compact) 
$3$-manifold 
that is a connected sum $(S^1\times S^2)\# M'$, 
and let $K\subset M$ be an unframed knot that crosses only once one of the spheres
${t}\times S^2\subset (S^1\times S^2)\# M'$, then $|K|=2$.
\end{thm}
\begin{emf}\label{proofExample}
Let $K_f$ be a framed knot corresponding to $K$. As it was essentially 
shown in the works of Hoste and Przytycki~\cite{HostePrzytycki} and later in the work of the author~\cite{ChernovLegendrian} Theorem 4.1.1, all the knots $K_f^{2i}, i\in \Z,$ are isotopic. Similarly all the knots $K_f^{2i+1}, i\in \Z,$ are also isotopic. 
Thus $|K|\leq 2$.

It is possible to show that $K_f^0$ and $K_f^1$ are not isotopic for all $K_f$ in all $M$, see for example~\ref{framedcomponents}, and thus $|K|\geq 2$. \qed
\end{emf}

\section{Some results needed for the proof of Theorem~\ref{aslk}}
\subsection{The covering $\pr:\mathcal F\rightarrow \mathcal
C$}\label{sectioncovering}

In this section the manifold $M$ is oriented but not necessarily compact.

Let $\mathcal C$ be a component of
the space of unframed curves in $M$. 
Let $\mathcal F'$ be a connected component of the space of framed curves in
$M$, such that the curves from $\mathcal F'$ realize curves from $\mathcal
C$ if we forget the framing.
Put $\pr':\mathcal F'\rightarrow
\mathcal C$ to be the forgetting of the framing mapping.

Let $p:\mathcal F'\to \mathcal F$ be the quotient by the 
following equivalence relation: $f'_1\sim f'_2$ if 
there exists a path $I:[0,1]\rightarrow \mathcal F'$ connecting $f'_1$ and
$f'_2$ such that $\Im(\pr'(I))=\pr'(f'_1)=\pr'(f'_2)$. (This means 
that we
identify two framed curves if the nonzero sections of the normal bundle to the curve induced by the framings are homotopic as nonzero sections.)
Put $\pr:\mathcal F\rightarrow \mathcal C$ to be the mapping such that
$\pr\circ p=\pr'$.

\begin{lem}\label{covering}
$\pr:\mathcal F\rightarrow \mathcal C$ is a covering with a structure group $\Z$, moreover this covering is normal, i.e.~$\Im(\pr_*(\pi_1(\mathcal F)))$
is a normal subgroup of $\pi_1(\mathcal C)$.

The  mapping $\delta:\pi_1(\mathcal C, c)\rightarrow \Z$, that maps the class $[\alpha]\in\pi_1(\mathcal C, c)$ of a loop $\alpha:[0,1]\to \mathcal C$ to
the element $\delta([\alpha])$ of the structure group $\Z$ of the covering such that $\delta([\alpha]) \cdot \widetilde \alpha(0)=\widetilde \alpha(1)$, is a homomorphism. (Here $\cdot$ denotes the action of the structure group $\Z$ and $\widetilde \alpha:[0,1]\to \mathcal F$ is a lift of $\alpha$.) 
\end{lem}

The proof of the Lemma is straightforward.

\begin{defin}[of isotopic knots from $\mathcal F$]\label{isotopy}
Let $K_0, K_1\in\mathcal F$ be such that
$\pr(K_0)$ and
$\pr(K_1)$ are knots (embedded curves). Then $K_0$ and $K_1$ are said to be
{\em isotopic\/} if there exists a path $q:[0,1]\rightarrow \mathcal F$ such that 
$q(0)=K_0, q(1)=K_1$, and $\pr\circ q$ is an isotopy (of unframed knots).
Lemma~\ref{covering} implies that framed knots $K_{f,0},  K_{f,1}\in\mathcal F'$
are framed isotopic if and only if $p(K_{f,0})$ and
$p(K_{f,1})$ are isotopic in $\mathcal F$.
\end{defin}

\begin{rem}\label{VassilievGoussarov}
Finite order invariants with values in an Abelian group $\mathcal A$ can be easily defined in the category of 
knots from $\mathcal F$. Clearly the mapping $p:\mathcal F'\to \mathcal F$ induces the natural isomorphism of the groups of order $\leq n$  Vassiliev-Goussarov invariants $p^*: V^n_{\mathcal F, \mathcal A}\to V^n_{\mathcal F', \mathcal A}$ of framed knots from $\mathcal F$ and of knots from $\mathcal F'$. 
\end{rem}

\subsection{$h$-principle and some facts about $\pi_1(\mathcal C)$}
\begin{emf}\label{h-principle}{\em $h$-principle for curves in $M$.\/}

Put $p:STM\rightarrow M$ to be the unit two-sphere tangent bundle of $M^3$.
The Smale-Hirsch $h$-principle, see for example~\cite{Gromov}, says that the space of curves in $M$ is weak
homotopy
equivalent to  the space of free (continuous) loops  $\Omega STM$ in $STM$. The weak
homotopy equivalence is given by mapping a curve $C$ to a loop $\vec
C\in\Omega STM$ that sends a point $t\in S^1$ to the point of $STM$
corresponding to the direction of the velocity vector of $C$ at $C(t)$.

A loop $\alpha\in\pi_1 (\Omega STM,\vec K)$ 
is a mapping $\mu_{\alpha}:T^2=S^1\times S^1\rightarrow STM$,
with $\mu_{\alpha}\big|_{1\times S^1}=\vec K$ and
$\mu_{\alpha}\big|_{S^1\times 1}$ being the
trace of 
$1\in S^1=\{z\in\C\big | |z|=1\}$ under the homotopy of $\vec K$ described by $\alpha$.
Put $t(\alpha)=\mu_{\alpha}\big|_{S^1\times 1}\in\pi_1(STM, \vec K(1))$.
Since $\pi_1(T^2)=\Z\oplus\Z$ is commutative, we get that
$t:\pi_1(\Omega STM,\vec K)\rightarrow \pi_1(STM, \vec K(1))$ is a
surjective
homomorphism of $\pi_1(\Omega STM,\vec K)$ onto the centralizer $Z(\vec K)$
of $\vec K\in\pi_1(STM,\vec K(1))$. ($t$ is surjective because if $\gamma\in
Z(\vec K)$, then we can take $\mu_{\alpha}:S^1\times S^1\rightarrow STM$ to
be such
that $\mu_{\alpha}\big|_{1 \times S^1}=\vec K$, $\mu_{\alpha}\big|_{
S^1\times 1}=\gamma$ and
$\mu_{\alpha}$ maps the $2$-cell of $S^1\times S^1$ according to the
commutation
relation between $\vec K$ and $\gamma$. Clearly $t(\alpha)=\gamma$ for
$\alpha\in\pi_1(\Omega STM, \vec K)$ that corresponds to $\mu_{\alpha}$.)

If for the loops 
$\alpha_1, \alpha_2\in\pi_1(\Omega STM, \vec K)$
we have $t(\alpha_1)=t(\alpha_2)\in\pi_1(STM, \vec K(1))$ , then the
mappings
$\mu_{\alpha_1}$ 
and $\mu_{\alpha_2}$ of $T^2$ corresponding to these loops can be deformed
to
be identical on the $1$-skeleton of $T^2$. Clearly the obstruction for
$\mu_{\alpha_1}$ and $\mu_{\alpha_2}$ to be homotopic as mappings of $T^2$
(with the mapping of the $1$-skeleton of $T^2$ fixed under homotopy)
is the element of $\pi_2(STM)$ obtained by gluing together the two $2$-cells
of the two tori along the common boundary. (In particular we get the
Proposition of
V.~L.~Hansen~\cite{Hansen}
that $t:\pi_1(\Omega X,\omega)\rightarrow Z(\omega)<\pi_1(X,\omega(1))$ is
an isomorphism, provided that $\pi_2(X)=0$.)

Since every oriented $3$-dimensional manifold is parallelizable, we can
identify $STM$ with $S^2\times M$ and $p:STM\rightarrow M$ with $p:S^2\times
M\rightarrow M$. Thus using the $h$-principle 
we can view $t$ as a surjective homomorphism
$t:\pi_1(\mathcal C, K)\rightarrow Z(K)<\pi_1(M, K(1))$.

Since $\pi_2(S^2\times M)=\Z\oplus\pi_2(M)$, we get that if $\alpha_1,
\alpha_2\in\pi_1(\mathcal C, K)$ are such that $t(\alpha_1)=t(\alpha_2)\in
Z(K)<\pi_1(M, K(1))$, then the obstruction for $\alpha_1$ and $\alpha_2$ to
be equal in $\pi_1(\mathcal C, K)$ is an element of $\Z\oplus\pi_2(M)$.
\end{emf}

\begin{prop}\label{centralizerfreeproduct}
Let $M$ be a connected sum $M_1\#M_2$, and let $K$ be a loop in $M$. 
\begin{description} 
\item[1]
If $K$ is not free homotopic to a loop $K'$ that is
contained either in $M_1\setminus B_3\subset M_1\#M_2$ or in $M_2\setminus
B_3\subset M_1\#M_2$, then the centralizer $Z(K)$ of $K\in\pi_1(M, K(1))$ is
an infinite cyclic group that contains $K\in\pi_1(M, K(1))$.
\item[2] If $K$ is a non-contractible 
loop that is contained in $M_1\setminus B_3\subset
M_1\#M_2$, then the centralizer $Z(K)$ of $K\in\pi_1(M_1\#M_2, K(1))$ is the
centralizer of $K\in\pi_1(M_1\setminus B_3, K(1))=\pi_1(M_1, K(1))$.
\end{description}
\end{prop}

To prove this Proposition one observes that by the Theorem of Van Kampen
$\pi_1(M_1\#M_2)$ is the free product of $\pi_1(M_1)$ and of $\pi_1(M_2)$.
After this the proof of the Proposition is an exercise in group theory.

The following Proposition is an immediate consequence of the Sphere Theorem,
see~\cite{Hempel}.

\begin{prop}\label{pi2}
Let $M$ be a (not necessarily closed) oriented manifold that is a connected sum 
$\# _{i=1}^k M_i$ of prime compact $3$-manifolds 
and let $\phi_i:S^2\rightarrow M$, $i=1,\dots, k,$ be
the embedded spheres with respect to which $M$ is the connected sum. Then
$\pi_2(M, x)$ as the $\pi_1$-module is generated by the classes of spheres 
$\phi_i$, i.e. every $s\in\pi_2(M,x)$
can be written as a product
$s=\prod _{i\in I} s_i^{\pm 1}$ 
of the spheroids $s_i$ such that 
\begin{description}
\item[a]
$s_i$ maps the the lower hemisphere of $S^2$ to a path connecting $x$ to
the south pole of $\phi_j(S^2)$, for some $j\in \{1, \dots, k\}$;
\item[b] $s_i$ maps the upper hemisphere of $S^2$ as it is described by
$\phi_j(S^2)$. (The upper hemisphere with all the points of the equator glued
together is naturally identified with $S^2$.)
\end{description}
\end{prop}

\begin{emf}\label{gamma2}{\em Loops $\gamma_1$, $\gamma_2$ and $\gamma_s$.\/}
Let $\mathcal C$ be a connected component of the space of curves in $M$ and
let $K\in\mathcal C$ be a knot. 

Let $\gamma_1$ be the isotopy of $K$ to
itself that is the sliding of $K$ along itself induced by the full 
rotation of the parameterizing circle.

Let $\gamma_2$ be the deformation of $K$ described in
Figure~\ref{obstruction.fig}. 

\begin{figure}[htbp]
 \begin{center}
  \epsfxsize 10cm
  \hepsffile{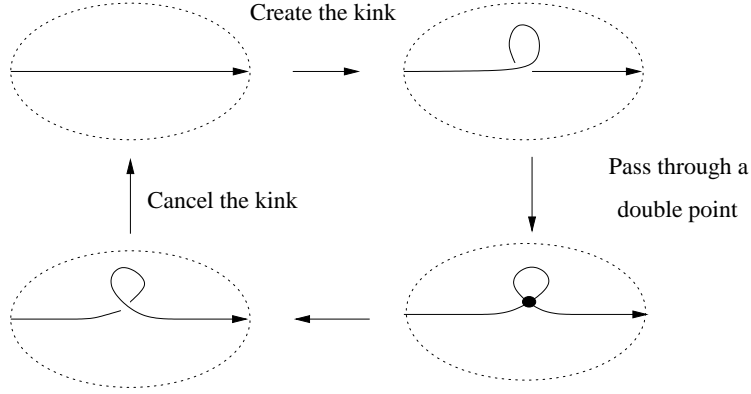}
 \end{center}
\caption{The loop $\gamma_2$.}
\label{obstruction.fig}
\end{figure}

Let $s\in\pi_2(M, K(1))$ be an element that is realizable by a mapping
$s:S^2\rightarrow M$ of the type described in Proposition~\ref{pi2}.
Put $\gamma_s\in\pi_1(\mathcal C, K)$ to be a loop, under which the knot $K$
does not move anywhere except of a small arc located close to $1\in S^1$. 
The points of the arc first slide along the boundary of the tubular neighborhood of
the path $\rho$. Then
the arc reaches the embedded sphere (that is $s$ restricted to the
upper hemisphere) and slides around the sphere, see Figure~\ref{gammas.fig}. 
Finally the points of the arc slide back along the boundary of the tubular 
neighborhood of the path $\rho$.

\begin{figure}[htbp]
 \begin{center}
  \epsfxsize 11cm
  \hepsffile{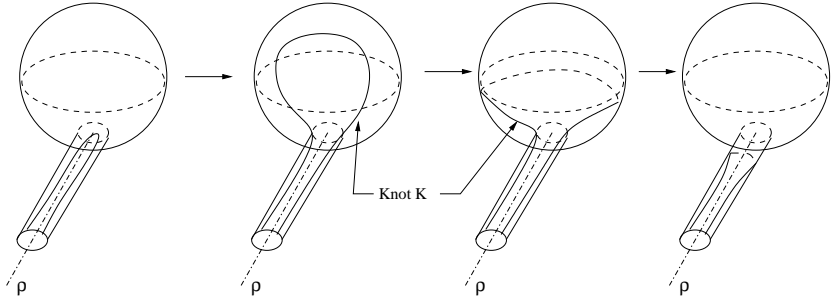}
 \end{center}
\caption{The loop $\gamma_s$.}
\label{gammas.fig}
\end{figure}

Let $t:\pi_1(\mathcal C, K)\rightarrow Z(K)<\pi_1(M, K(1))$ be the
homomorphism described in~\ref{h-principle}. 
Let $\alpha_1, \alpha_2\in\pi_1(\mathcal C, K)$
be such that $t(\alpha_1)=t(\alpha_2)\in\pi_1(M, K(1))$. Let
$m\oplus\epsilon\in\Z\oplus\pi_2(M)=\pi_2(S^2\times M)=\pi_2(STM)$ 
be the obstruction for $\alpha_1$ and $\alpha_2$ to be equal in
$\pi_1(\mathcal C, K)$,
see~\ref{h-principle}. The Sphere Theorem~\cite{Hempel} says 
that $\epsilon \in\pi_2(M)$
can be realized as a product $\epsilon=\prod_{i=1}^k s_i^{\sigma i}$ 
of the spheroids of the type described above ($\sigma _i=\pm 1$). 
Straightforward geometric considerations show that $t(\alpha_1 \prod_{i=1}^k
\gamma_{s_i}^{\sigma i})=t(\alpha_2)$ and the obstruction for $\alpha_1
\prod_{i=1}^k \gamma_{s_i}^{\sigma i}$ and $\alpha_2$ to be equal
elements in $\pi_1(\mathcal C, K)$ is an element $m'\oplus
0\in\Z\oplus\pi_2(M)=\pi_2(STM)$. Finally one verifies that 
$t(\alpha_1 (\prod_{i=1}^k \gamma_{s_i}^{\sigma i})\gamma_2^{m'})=t(\alpha_2)$
and the obstruction for $\alpha_1 (\prod_{i=1}^k \gamma_{s_i}^{\sigma
i})\gamma_2^{m'}$ and $\alpha_2$ to be homotopic vanishes, i.e. they realize
the same elements of $\pi_1(\mathcal C, K)$.

Thus 
we get the following Lemma.
\end{emf}

\begin{lem}\label{obstructions}
Let $M$ be a (not necessarily closed) oriented $3$-manifold that is not $(S^1\times S^2)\#M'$, 
let $\mathcal C$ 
be a connected
component of the space of curves in $M$, and let $K\in\mathcal C$ be a knot.
Let $\alpha_1, \alpha_2\in\pi_1(\mathcal C,K)$ be such that
$t(\alpha_1)=t(\alpha_2)\in\pi_1(M, K(1))$. 
\begin{description}
\item[1]
Then there exist 
spheroids $s_i\in\pi_2(M, K(1))$, $i=1, \dots k$, 
of the type described above in~\ref{pi2}, $\sigma_i=\pm 1$, and $m\in\Z$, such that 
$\alpha_1(\prod_{i=1}^k \gamma_{s_i}^{\sigma
i})\gamma_2^{m}=\alpha_2\in\pi_1(\mathcal C, K)$.
\item[2] If $\pi_2(M)=0$, then $\alpha_1\gamma_2^m=\alpha_2$, for some
$m\in\Z$.
\end{description}

\end{lem}

\section{Proof of Theorem~\ref{aslk}}\label{Proofaslk}
Since both $\aslk$ and $\wt \aslk$ are clearly finite order invariants, it suffices to show, see~\ref{VassilievGoussarov}, that the corresponding finite order invariants of knots from $\mathcal F$ are well-defined.

Invariants $\wt \aslk$ and $\aslk$ are constructed in a similar fashion. 
To construct $\wt \aslk$ choose a preferred knot $K_f\in\mathcal F$ and a value 
$\wt \aslk (K_f)$. (The choice of  $\wt \aslk (K_f)$ is the only ambiguity in the construction of $\wt \aslk$ and it corresponds to the fact that $\wt \aslk$ is uniquely defined by its properties up to an additive constant.)

Let $\overline K_f\in \mathcal F$ be a knot, and let $\rho:[0,1]\to \mathcal F$ be a generic path such that $\rho(0)=K_f$ and $\rho(1)=\overline K_f$. Put $J_{\rho}\subset [0,1]$ to be set of instances when the knot becomes singular under the deformation $\rho$. At these instances $\rho$ crosses the
discriminant in $\mathcal F$. (The {\em discriminant\/} is the
subspace of $\mathcal C$ formed by singular knots.) Since $\rho$ is
generic, at these instances the knot has one transverse double point, that
separates the knot into two loops. Put $\wt J_{\rho}\subset 
J_{\rho}$ to be
set of instances when one of the two loops of the singular knot is
contractible. Put $\sigma _j$, $j\in J_{\rho}$, to be the signs of the corresponding crossings of the discriminant.
(The sign of the crossing is the sign of the resolution of the
double point that occurs during the crossing.)
Put 
\begin{equation}\label{defindeltatildeaslk}
\Delta_{\wt \aslk}(\rho)=\sum_{j\in \wt J_{\rho}}2\sigma _{j}
\end{equation} and 
put 
\begin{equation}\label{defindeltaaslk}
\Delta_{\aslk}(\rho)=\sum_{j\in J_{\rho}}2\sigma _{j}.
\end{equation}

Clearly if the $\wt \aslk$ invariant of Theorem~\ref{aslk} invariant does exist then $\wt \aslk (\overline K_f)-\wt \aslk (K_f)=\Delta_{\wt \aslk}(\rho)$ and we put $\wt \aslk (\overline K_f)=\wt \aslk(K_f)+\Delta_{\wt \aslk}(\rho)$.

To show that  $\wt \aslk$ is well-defined we have to prove that $\Delta_{\wt \aslk}(\rho)=\Delta_{\wt \aslk}(\bar \rho)$ for any other generic path $\bar \rho:[0,1]\to \mathcal F$ with $\bar \rho (0)=K_f$ and $\bar \rho (1)=\overline K_f$. Since $\Delta_{\wt \aslk}(\rho^{-1})=-\Delta_{\wt \aslk}(\rho)$, it follows that to prove the existence of $\wt \aslk$ it suffices to show that $\Delta_{\wt \aslk}(\bar \alpha)=0$, for every generic closed loop $\bar \alpha:[0,1]\to \mathcal F$ with $\bar \alpha(0)=\bar \alpha(1)=K_f$. 

Similar considerations show that to prove the existence of the $\aslk$ invariant it suffices to show that $\Delta_{\aslk}(\bar \alpha)=0$, for every generic closed loop $\bar \alpha:[0,1]\to \mathcal F$ with $\bar \alpha(0)=\bar \alpha(1)=K_f$.

The codimension two (with respect to $\mathcal F$) stratum of the discriminant consists of singular knots with two distinct transverse double points. It is easy to see that
$\Delta_{\wt \aslk}(\alpha')=\Delta_{\wt  \aslk}(\alpha')=0$ and $\Delta_{\aslk}(\alpha')=\Delta_{\aslk}(\alpha')=0$,
for every small generic loop $\alpha'$ going around the
codimension two stratum. This implies (cf. Arnold~\cite{Arnoldcurves}) 
that if $\gamma$ is a generic loop in $\mathcal F$ that starts
at a nonsingular knot $K_f$,
then $\Delta_{\wt \aslk}(\bar \alpha)$ and $\Delta_{\aslk }(\bar \alpha)$ 
depend only on the element of 
$\pi_1(\mathcal F, K_f)$ realized by a generic loop $\bar \alpha$. 
Now it is clear that $\Delta_{\wt \aslk}:\pi_1(\mathcal F, K_f)\rightarrow \Z$ that maps the class of a generic loop $\bar \alpha$ (starting at $K_f$) to $\Delta_{\wt \aslk}(\bar \alpha)$ is a
homomorphism. Similarly $\Delta_{\aslk}:\pi_1(\mathcal F, K_f)\rightarrow \Z$ is also a homomorphism. 
To prove the Theorem it suffices to show that the homomorphisms $\Delta_{\wt \aslk}, \Delta_{\wt \aslk}: \pi_1(\mathcal F, K_f)\rightarrow \Z$ are zero-homomorphisms (under the corresponding conditions on the component $\mathcal F'$).

Consider the covering $\pr :\mathcal F\to \mathcal C$, see~\ref{covering}, and 
put $K=\pr (K_f)$. Similarly to the above we introduce homomorphisms $\Delta_{\wt \aslk}, \Delta_{\aslk}: \pi_1(\mathcal C, K)\rightarrow \Z$. (We will use $\Delta_{\wt \aslk}$ and $\Delta_{ \aslk}$ as the notation for homomorphisms from both $\pi_1(\mathcal F, K_f)$ and $\pi_1(\mathcal C, K)$.)  Clearly $\Delta_{\wt \aslk}(\bar \alpha)=\Delta_{\wt \aslk}(\pr(\bar \alpha))$ and $\Delta_{\aslk}(\bar \alpha)=\Delta_{ \aslk}(\pr(\bar \alpha))$, for all $\bar \alpha\in \pi_1(\mathcal F)$. 
A loop $\alpha\in\pi_1(\mathcal C)$ is liftable to $\mathcal F$ if and only if $\delta(\alpha)=0$, see~\ref{covering}.

Now the proof of Theorem~\ref{aslk} is finished modulo the following technical Theorem (with a rather hard proof).

\begin{thm}\label{technicaltheorem}
Let $M$ be closed not realizable as $M'\# (S^1\times S^2)$.
\begin{description}
\item[1] If $K$ is not homotopic to a loop contained inside of one of the irreducible summands of $M=\#_{j\in J}M_j$, then  
$\delta$ and $\Delta_{\aslk}:\pi_1(\mathcal C, K)\to \Z$ are equal homomorphisms.

\item[2] If $K$ is homotopic to a loop contained inside of one of the irreducible summands of $M=\#_{j\in J}M_j$, then 
$\delta$ and $\Delta_{\wt \aslk}:\pi_1(\mathcal C, K)\to \Z$ are equal homomorphisms.
\end{description}
\end{thm}

\subsection{Proof of Theorem~\ref{technicaltheorem}}

Observe that in the Proof of Theorem~\ref{technicaltheorem}  we can assume that $K$ is any nonsingular knot from  $\mathcal C$ we like.
Since $\Z$ is torsion free, to prove Statement {\bf 1} it suffices to show that for every $\alpha\in\pi_1(\mathcal C, K)$, there exists $i\neq 0$ such that $\Delta_{\wt \aslk}(\alpha^i)=\delta(\alpha^i)$. A similar observation 
holds for the proof of statement {\bf 2}.

We will need the following Lemma.

\begin{lem}\label{betaforobstructions}

Let $M$ be (a not necessarily closed) 
oriented manifold that is not $(S^1\times S^2)\#M'$.
Let $\gamma_1,\gamma_2, \gamma_s$ be the loops described in~\ref{gamma2}. Then 
\begin{description}
\item[1]
$\delta(\gamma_1)=\Delta_{\wt\aslk}(\gamma_1)=\Delta_{\aslk}(\gamma_1)=0$, and
$\delta(\gamma_2)=\Delta_{\wt\aslk}(\gamma_2)=\Delta_{\aslk}(\gamma_2)=2$ (see~\ref{covering},~\eqref{defindeltatildeaslk},~\eqref{defindeltaaslk}  for the definitions of $\delta, \Delta_{\wt \aslk}, \bar\Delta_{\aslk}$); 
\item[2]
$\delta(\gamma_s)=\Delta_{\aslk}(\gamma_s)=0$; 
\item[3]Let $M\neq (S^1\times S^2)\#M'$ be closed and oriented.
Let $K$ be a knot that is contained in one the irreducible summands of $M$,
then $\Delta_{\wt\aslk}(\gamma_s)=0$.
\end{description}
\end{lem}

\begin{emf}{\em Proof of Lemma~\ref{betaforobstructions}.\/}
The proofs of statement {\bf 1} and of the identity $\delta(\gamma_s)=0$
are obtained by straightforward geometric considerations.

Let $s'$ be the embedded sphere used to construct $\gamma_s$.
To get $\Delta_{\aslk}(\gamma_s)=0$ we observe that $\Delta_{\aslk}(\gamma_s)$ is
equal to the intersection index of $K\in H_1(M)$ and $s'\in H_2(M)$. (The moments
when the knot becomes singular are those when the branch of $K$ that slides
around $s'$ passes through a branch of $K$ that intersect $s'$.) On the
other hand, since $M\neq (S^1\times S^2)\#M'$, every embedded sphere separates
$M$ into two disjoint parts, and thus the intersection index of $K$ and $s'$ is
zero.

To get statement {\bf 3} we observe that since $K$ lies in the irreducible
summand of $M$, the only crossings of the discriminant that occur under
$\gamma_s$ are those when the branch that slides around the sphere $s'$
passes through the small neighborhoods of the crossings of the path $\rho$
(that was used to construct $\gamma_s$) and the sphere $s'$. Each small
neighborhood of the crossing of $\rho$ and $s'$ gives rise to two crossings
of the discriminant under $\gamma_s$. A straightforward verification show
that the signs of the two crossings are opposite, and that if one appearing
singular knot in such a pair has a contractible loop, then so does the other
one. Thus $\Delta_{\wt\aslk}(\gamma_s)=0$.\qed
\end{emf}

\begin{emf}
{\bf Let us prove statement {\bf 1} of Theorem~\ref{technicaltheorem}.}
Let $t:\pi_1(\mathcal C, K)\rightarrow
Z(K)<\pi_1(M, K(1))$ be the
homomorphism described in~\ref{h-principle}. 
Clearly $t(\gamma_1)=K\in\pi_1(M, K(1))$ for the loop $\gamma_1$ introduced
in~\ref{gamma2}.  
Take $\alpha\in\pi_1(\mathcal C, K)$. 
Proposition~\ref{centralizerfreeproduct}.1 implies
that that the centralizer $Z(K)$ of $K\neq 1\in\pi_1(M, K(1))$ is an
infinite cyclic group. 
Thus there exists a nonzero $i\in\Z$ such that $t(\alpha^i)=t(\gamma_1^k)$,
for some $k\in\Z$. By Lemma~\ref{obstructions} $\alpha^i=\gamma_1^k\prod
\gamma_{s_i}^{\sigma_i}\gamma_2^m$, for some $m\in\Z$, $\sigma_i=\pm 1$, and spheroids $s_i$ of the type described in~\ref{gamma2}.

By Lemma~\ref{betaforobstructions} we have 
$\delta(\gamma_1^i)=\Delta_{\wt\aslk}(\gamma_1^i)=0$, 
$\delta(\prod \gamma_{s_i}^{\sigma_i})=\Delta_{\wt\aslk}(\prod
\gamma_{s_i}^{\sigma_i})=0$, and $\delta(\gamma_2^m)=\Delta_{\wt\aslk}(\gamma_2^m)=2m$.
This implies that  
$\delta(\alpha^i)=\Delta_{\wt\aslk}(\alpha^i)$ for our choice of $i\neq 0$ and this finishes the proof of
statement {\bf 1}.

{\bf Below we prove statement {\bf 2}.} {\em Consider first the case of contractible $K$.\/} Without the loss of generality we can assume that $K$ is the unknot contained in a small ball $B\subset M$. Take $\alpha\in \pi_1(\mathcal C, K)$. Let $\overline \alpha \in\pi_1(\mathcal C, K)$ be the deformation of $K$ induced by the ambient diffeotopy of $M$ that is the sliding of $B$ (with $K$ in it) along the loop $t(\alpha)$. Clearly $\delta(\overline\alpha)=\Delta_{\wt\aslk}(\overline \alpha)=0$ and $\alpha=\overline \alpha \prod
\gamma_{s_i}^{\sigma_i}\gamma_2^m$, for some $m\in\Z$, $\sigma_i=\pm 1$, and spheroids $s_i$ of the type described in~\ref{gamma2}. After this the proof follows immediately similarly to the case of Statement {\bf 1}.

{\em Below in the Proof of the Theorem we assume that $K$ is not contractible.}

Let $M_1$ be the irreducible summand of $M$ that contains $K$.  
Proposition~\ref{centralizerfreeproduct}.2 implies that the centralizer of
$K\in\pi_1(M, K(1))$ is canonically isomorphic to the centralizer of
$K\in\pi_1(M_1, K(1))$. Thus the $h$-principle, see~\ref{h-principle},
implies that for any $\alpha\in \pi_1(\mathcal C, K)$ 
there exists 
$\bar \alpha \in\pi_1(\mathcal C, K)$ with
$t(\alpha)=t(\bar \alpha)\in\pi_1(M, K(1))$ 
such that $\bar \alpha$ is realizable
by a deformation of $K$ (in the space curves) 
that does not take $K$ outside of $M_1$.

Lemma~\ref{obstructions} says that $\alpha=\bar\alpha(\prod
\gamma_{s_i}^{\sigma_i})\gamma_2^m$, for some $m\in\Z, \sigma_i=\pm 1$, and
spheroids $s$ of the type described in~\ref{gamma2}. 

By Lemma~\ref{betaforobstructions} we get that 
$\delta(\alpha^i)=2mi +\delta(\bar\alpha^i)$ and $\Delta_{\wt\aslk}(\alpha^i)=2mi
+\Delta_{\wt\aslk}(\bar\alpha^i)$, for any nonzero $i\in\Z$.

Now statement {\bf 2} follows immediately from the
following Theorem~\ref{irreducible}. (This finishes the Proof of
Theorem~\ref{technicaltheorem} and of Theorem~\ref{aslk} modulo the Proof of
Theorem~\ref{irreducible}.)\qed
\end{emf}

\begin{rem}\label{aslkcompactnotclosed}
Every compact manifold also admits a decomposition into a connected sum of prime manifolds, see for example~\cite{Jaco} exercise {\textrm II}.18. 
Thus the above proof of statement {\bf 1} of Theorem~\ref{aslk} holds for compact (not necessarily closed) $M\neq (S^1\times S^2)\#M'$ decomposed as a sum of primes.
\end{rem}

\begin{thm}\label{irreducible}
Let $M$ be a closed oriented irreducible $3$-manifold, let $\mathcal C$
be a connected component of the space of curves in $M$ that consists of not
contractible curves, and let
$K\in\mathcal C$ be a knot. Then for any $\alpha\in\pi_1(\mathcal C, K)$
there exists a nonzero $i\in\Z$ such that
$\delta(\alpha^i)=\Delta_{\wt\aslk}(\alpha^i)$.
\end{thm}

\subsection{Proof of Theorem~\ref{irreducible}}\label{proofirreducible} 
By the Sphere Theorem~\cite{Hempel} $\pi_2(M)=0$.
Let $t:\pi_1(\mathcal C, K)\rightarrow Z(K)<\pi_1(M, K(1))$ be the
homomorphism introduced in~\ref{h-principle}.

If $\pi_1(M)$ is finite, then for any $\alpha\in\pi_1(\mathcal C, K)$ there
exists $i\neq 0$ such that $t(\alpha^i)=1=t(1)$. Lemma~\ref{obstructions}.2
says that there exists $m\in\Z$ such that $\alpha^i=\gamma_2^m$. Now
Lemma~\ref{betaforobstructions} implies that
$\delta(\alpha^i)=2m=\Delta_{\wt\aslk}(\alpha^i)$.
 This finishes the proof in the case of $\pi_1(M)$ being finite, and {\em below in
the proof we assume that $\pi_1(M)$ is infinite.\/} 

Since $\pi_1(M)$ is
infinite and $M$ is orientable and irreducible, the result of
Epstein, see~\cite{Epstein} or~\cite{Hempel} Corollary 9.9, implies
{\em that $\pi_1(M)$
is torsion free.\/} 

Let $p:STM=S^2\times M\rightarrow M$ be the $S^2$-fibration.
To a loop $\alpha$ in $\mathcal C$ such that $\alpha(1)=K$ we correspond a
mapping $\mu_{\alpha}:S^1\times S^1\rightarrow STM$ as it is described
in~\ref{h-principle}. Put $\bar\mu_{\alpha}=p(\mu_{\alpha}):S^1\times S^1\rightarrow
M$.

\begin{defin} Let $M'$ be a compact oriented 
$3$ manifold and let $F\neq S^2$ be a surface. A map
$\Psi:F\rightarrow M'$ is {\em essential\/} if 
$\ker(\Psi_*:\pi_1(F)\rightarrow\pi_1(M'))=1$.
\end{defin}

\begin{lem}\label{nonessential}
Let $M'$ be a submanifold of an oriented closed irreducible $M$ 
with infinite $\pi_1(M)$ such that $\pi_2(M')=0$ and 
the inclusion $i:M'\rightarrow M$
induces the injective homomorphism $i_*:\pi_1(M')\rightarrow \pi_1(M)$. 
Let $\alpha$ be a loop in the space of curves in $M$ such that $\Im \bar
\mu_{\alpha}\subset M'$ and 
$\bar \mu_{\alpha}:S^1\times S^1\rightarrow M'$ is not essential.
Then $\Delta_{\wt\aslk}(\alpha)=\delta(\alpha)$ (for the homomorphisms $\Delta_{\wt\aslk}$ and
$\delta$ defined with respect to the ambient manifold $M$).
\end{lem}

\begin{emf}{\em Proof of Lemma~\ref{nonessential}.\/}
Since $\ker\bar\mu_{\alpha_*}\neq 1$ and $\pi_1(S^1\times S^1)=\Z\oplus\Z$ is
generated by $1\times S^1$ and $S^1\times 1$, we have that
$K^j=(\bar\mu_{\alpha_*}(1\times S^1))^j=(\bar\mu_{\alpha_*}(S^1\times
1))^i=t(\alpha)^i\in\pi_1(M', K(1))$, for
some $i,j\in\Z$ (that are not both zero). Since $\pi_1(M)$ is torsion free and
$i_*:\pi_1(M')\rightarrow \pi_1(M)$ is 
injective, we have that $\pi_1(M')$ is torsion free.
Thus both $i$ and $j$ are nonzero.

Clearly $t(\gamma_1)=K\in\pi_1(M,
K(1))$ for the loop $\gamma_1$ introduced in~\ref{gamma2}.  
Since $\pi_2(M')=0$, Lemma~\ref{obstructions}.2
implies that there exists $m\in\Z$ such that $\alpha^i=\gamma_1^j\gamma_2^m$.
Lemma~\ref{betaforobstructions} says that
$\Delta_{\wt\aslk}(\gamma_1)=\delta(\gamma_1)=0$ and $\Delta_{\wt\aslk}(\gamma_2)=\delta(\gamma_2)=2$.
Thus $\Delta_{\wt\aslk}(\alpha^i)=\delta(\alpha^i)$ and since $i\neq0$ we get the
statement of the Lemma.\qed 

\end{emf}

\begin{rem}\label{nonessentialaslk}
Lemma~\ref{nonessential} proves the Theorem for nonessential 
$\bar \mu_{\alpha}$.
{\em Below in the Proof we assume that $\bar\mu_{\alpha}$ is an essential
mapping.\/}

Clearly the proof above can be easily modified to show that if the 
prime (not necessarily closed) summand of compact $M\neq (S^1\times S^2)\#M'$ that contains $K\in \mathcal C$ does not have essential tori, then both $\wt \aslk$ and $\aslk$ invariants are defined for all knots in $\mathcal F'$.
\end{rem}

\begin{defin} A surface $F\neq S^2$, properly embedded in a $3$-manifold
$M'$ (or embedded into $\p M'$), is {\em compressible,\/} if there exists a
disc $D^2\subset M$ such that $D^2\cap F=\p D^2$ and $\p D^2$ is not homotopically
trivial in $F$; otherwise $F$ is {\em incompressible\/} in $M$. A compact
orientable irreducible $3$-manifold is {\em Haken\/} if it contains a
two-sided incompressible surface.
\end{defin}

\begin{defin} Let $(\mu, \nu)$ be a pair of relatively prime integers. Let
$D^2=\{(r, \theta); 0\leq r\leq 1, 0\leq \theta\leq 2\pi\}$ be the unit disc
with polar coordinates. A {\em fibered torus of type $(\mu, \nu)$\/} is the
quotient of $D^2\times[0,1]$ via
$((r,\theta),1)=((r,\theta+\frac{2\pi\nu}{\mu}),0)$. The fibers are the
closed curves that are the unions $((r,\theta)\times
[0,1])\cup((r,\theta+\frac{2\pi\nu}{\mu})\times
[0,1])\cup((r,\theta+\frac{4\pi\nu}{\mu})\times [0,1])\cup\dots$ etc. for
fixed $(r,\theta)$. If $|\mu|>1$, then the solid torus is said to be
{\em exceptionally fibered\/}, the core of the torus is the {\em exceptional
fiber,\/} and $\mu$ is the {\em index\/} of the exceptional fiber. Otherwise the torus is {\em regularly fibered} and each fiber is a
{\em regular fiber.\/}

An orientable $3$-manifold $M'$ is called {\em Seifert-fibered\/} if it is a
union of pairwise disjoint closed curves (fibers), such that each fiber
has a neighborhood consisting of fibers that is homeomorphic to a fibered
solid torus via a fiber preserving homeomorphism.  A fiber of a Seifert
manifold $M'$ is exceptional if its neighborhood is homeomorphic to an
exceptionally fibered solid torus via a fiber preserving homeomorphism. 

The quotient space obtained from $M'$ via identifying all points in the
fiber for all the fibers is the {\em orbit space\/} and the images of the
exceptional fibers are the {\em cone points.\/}

\end{defin}

\begin{defin}[of a characteristic submanifold]\label{characteristic} 
A codimension zero submanifold $S$ of a closed oriented manifold
$M$ is called a {\em characteristic submanifold\/} if 
\begin{description}
\item[1]
each component $X$ of
$S$ admits a structure of a Seifert-fibered space or of a total space of a
$[0,1]$-bundle;
\item[2] if $W$ is a nonempty codimension zero submanifold of $M$ that
consists of components of $\overline{M\setminus S}$, then $S\cup W$ does not
satisfy $1$;
\item[3] if $S'$ is a codimension zero submanifold of $M$ that satisfies $1$
and $2$, then $S'$ can be deformed into $S$ by a proper isotopy.
\end{description}
\end{defin}

\begin{rem}\label{torussubstitution} 
Since $\delta, \Delta_{\aslk}$ and $\Delta_{\wt\aslk}$ are homomorphisms to an abelian group $\Z$,  they can be considered as homomorphisms $H_1(\mathcal C)\to \Z$. In particular the values of the homomorphism depend only on the conjugacy class of a loop in $\pi_1(\mathcal C)$ and  further in the proof we are free to change the mapping $\mu_{\alpha}:S^1\times S^1\to STM$ to 
any mapping homotopic to it.
Using this observation and Lemmas~\ref{obstructions}  and~\ref{betaforobstructions} we get that in fact further in the proof we are free to substitute instead of the mapping $\bar \mu_{\alpha}=p(\mu_{\alpha}):S^1\times S^1\to M$ any mapping of the torus homotopic to it.
\end{rem}

\begin{emf}\label{essential}{\em Let us reduce the proof of the Theorem in
the case of essential $\bar\mu_{\alpha}$, to the
case where the homotopy $\alpha$ happens
inside a Seifert-fibered submanifold $S$ of $M$ with $\pi_2(S)=0$.}

The Torus Theorem by Casson-Jungreis~\cite{CassonJungreis} and
Gabai~\cite{Gabaitorus} says that since $\bar\mu_{\alpha}$ is essential,
then either $M$ contains an embedded incompressible torus, and thus is Haken,
or it is a Seifert-fibered space. If it is Seifert-fibered then we have
made the reduction, so consider the case where $M$ is Haken.
%

The results of Jaco and Shalen~\cite{JacoShalen} and
Johannson~\cite{Johannson}, Proposition 9.4, say that 
$M$ has a well-defined characteristic submanifold $S$. 
The Enclosing Theorem by Jaco and Shalen~\cite{JacoShalen} and
Johannson~\cite{Johannson}, Theorem 12.5, says that
since $M$ is Haken there exists a mapping $\lambda:S^1\times S^1\rightarrow
S\subset M$ such that 
$\lambda$ is homotopic to $\bar\mu_{\alpha}$. 
There are only two oriented total spaces of $[0,1]$-bundles over surfaces that admit
essential mappings of $S^1\times S^1$. 
They are the $[0,1]\times S^1\times S^1$
and the unique $[0,1]$-bundle over the Klein-bottle with oriented total space.
Both these two spaces admit the structure of a Seifert-fibered space. (In
the case of a $[0,1]$-bundle over the Klein bottle it is a Seifert-fibered
space over a Moebius strip.) Thus we get that $\lambda$ is a map into a
Seifert-fibered component of $S$. 

Now Remark~\ref{torussubstitution} and $h$-principle~\ref{h-principle} allow us to assume that $\bar \mu_{\alpha}$ is contained in the Seifert-fibered component of $S$.

%

Finally to complete the reduction we observe that clearly the value of
$\delta(\alpha)$ does not depend on whether we regard $\delta$ as a
homomorphism from $\pi_1$ of $\mathcal C$ or of $\mathcal C_S$. (Here $\mathcal C_S$ is the component of the space of all curves in $S$ that contains $K$.) The inclusion of every 
component of $S$ into $M$ induces a monomorphism of fundamental groups,
(see~\cite{Johannson} remark on p.27 and 8.2). Thus the value of
$\Delta_{\wt\aslk}(\alpha)$ also does not depend on whether we regard $\Delta_{\wt\aslk}$ as a
homomorphism $\Delta_{\wt\aslk}:\pi_1(\mathcal C, K)\rightarrow \Z$ or as 
$\Delta_{\wt\aslk}:\pi_1(\mathcal C_S, K)\rightarrow \Z$. Also since the
homotopy $\alpha'$ happens inside a Seifert-fibered component of the
characteristic submanifold $S$ of irreducible $M$, we
have that $\pi_2(S)=0$. 
\end{emf}

{\em Thus we have reduced the proof of the Theorem to the case where $M$ is an
oriented connected compact (not necessarily closed) Seifert-fibered manifold,
with $\pi_2(M)=0$, $\pi_1(M)$ is torsion free,  
$\bar \mu_{\alpha}$ is an essential torus in $M$ (see
Lemma~\ref{nonessential}),
and $K$ is a non-contractible knot in $M$.\/} 

\begin{rem}\label{wtaslkcompacthaken}
A straightforward verification shows that the proof of Theorem~\ref{irreducible} holds under the assumption that $M$ is a compact (rather than closed) prime 
manifold that either does not admit essential tori; or is Haken, so that $\bar \mu_{\alpha}$ is homotopic into a characteristic submanifold.   
\end{rem}

\begin{defin} Let $q:S\rightarrow F$ be a Seifert-fibration. A mapping
$\lambda:S^1\times S^1\rightarrow S$ is said to be {\em vertical\/} 
with respect to $q$ if $q^{-1}(q\lambda(S^1\times S^1))
=\lambda (S^1\times S^1)$ and
$\lambda (S^1\times S^1)$ does not contain exceptional fibers of
$q:S\rightarrow M$.
\end{defin}

{\em The following three Lemmas 
prove the Theorem in the case where $\bar\mu_{\alpha}$
is homotopic to a vertical mapping of $S^1\times S^1$.\/}

\begin{lem}\label{commute}
Let $p:S\rightarrow Y$ be a locally trivial $S^1$-fibration of an oriented
manifold $S$ over a (not necessarily orientable) manifold $Y$. Let
$f\in\pi_1(S)$ be the class of an oriented
$S^1$-fiber of $p$, and let $\alpha$ be an
element of $\pi_1(S)$. Then:
\begin{description}
\item[a] $\alpha f=f\alpha\in\pi_1(S)$, provided that $p(\alpha)$ is an
orientation preserving loop in $Y$;
\item[b] $\alpha f=f^{-1}\alpha\in\pi_1(S)$, provided that $p(\alpha)$ is an
orientation reversing loop in $Y$.
\end{description}
\end{lem}

\begin{emf}{\em Proof of Lemma~\ref{commute}.\/}
If we move an oriented fiber along the loop $\alpha\in S$, then in the end
it comes to itself either with the same or with the opposite orientation.
It is easy to see that it comes to itself with the opposite orientation if
and only if $p(\alpha)$ is an orientation reversing loop in $Y$.\qed
\end{emf}

\begin{lem}\label{locallytrivial}
Let $\tilde p:N\rightarrow G$ be a Seifert-fibration, let $F\subset G$ be a
connected submanifold with boundary that does not contain cone points. Let
$M=\tilde p^{-1}(F)\subset N$, and let $p:M\rightarrow F$ be the
corresponding locally trivial $S^1$-fibration. Let $\alpha$ be a homotopy of
a non-contractible knot $K$ such that $\bar \mu_{\alpha}\subset M$, then
$\Delta_{\wt\aslk}(\alpha)=\delta(\alpha)$.
\end{lem}

\begin{emf}\label{prooflocallytrivial}{\em Proof of Lemma~\ref{locallytrivial}.\/}
Put $f$ to be the class of
the oriented fiber of $p:M\rightarrow F$.
Let $t:\pi_1(\mathcal C, K)\rightarrow Z(K)<\pi_1(M, K(1))$ be the 
homomorphism introduced
in~\ref{h-principle}.

{\bf Consider the case of $p(K)\neq 1\in\pi_1(F)$.\/} 
Clearly $p_*t(\alpha)\in Z(p_*(K))<\pi_1(F)$. Since $\p F\neq\emptyset$,
we have that $\pi_1(F)$ is a free group. Since $p(K)\neq 1\in\pi_1(F)$, we
get that there exists $i\neq 0$ and $j$ such that
$p_*t(\alpha^i)=p_*(K^j)\in\pi_1(F)$. Using Lemma~\ref{commute} and the fact
that $f$ generates $\ker p_*$, we obtain that there exists $k$ such that
\begin{equation}\label{express}
t(\alpha^i)=K^jf^k\in\pi_1(M, K(1)).
\end{equation} 

{\em Consider the case of $p(K)$ being an orientation reversing loop on
$F$.\/} Since $t(\alpha^i)$ commutes with $K$ in $\pi_1(M)$ and by
Lemma~\ref{commute} $fK=Kf^{-1}$, we get that $f^{2k}=1$. Since
$\pi_2(M)=0$, we have that $f$ has infinite order.
Thus $k=0$ and $t(\alpha^i)=K^j$. 

Let $\gamma_1$ be the isotopy of $K$ to itself introduced in~\ref{gamma2}.
Clearly $t(\gamma_1)=K$. Thus by Lemma~\ref{obstructions}
we get that $\alpha^i=\gamma_1^j\gamma_2^s$, for some $s\in\Z$. Using
Lemma~\ref{betaforobstructions}
we get that
$\Delta_{\wt\aslk}(\alpha^i)=\delta(\alpha^i)$. Since $i$ was chosen to be nonzero we
get that $\Delta_{\wt\aslk}(\alpha)=\delta(\alpha)$.

{\em Consider the case of $p(K)$ being an orientation preserving loop on
$F$.\/} Since $M$ is orientable, we get that the $S^1$-fibration over $S^1$
(parameterizing the knots) induced from $p$ by $p\circ K:S^1\rightarrow F$ is
trivializable. Hence we can coherently orient the fibers of the induced 
fibration $\bar p:S^1\times S^1\rightarrow S^1$. The orientation of the $S^1$-fiber 
of $\bar p$ over $t\in S^1$ induces the orientation of the $S^1$-fiber of $p$ that contains 
$K(t)$. Let $\gamma_3$ be the homotopy of $K$ that slides every 
point $K(t)$ of $K$ inside the fiber that contains $K(t)$ with unit velocity
in the direction specified by the
orientation of the fiber of $\bar p$ over $t\in S^1$.
Clearly $t(\gamma_3)=f$. Thus $\alpha^i=\gamma_1^j\gamma_3^k\gamma_2^s$, for
some $s\in\Z$. 

Let us show that $\Delta_{\wt\aslk}(\gamma_3)=0$. The only singular knots
that arise under $\gamma_3$ are those that have a double point projecting
to a double point $d$ of $p(K)$. Every double point $d$ of $p(K)$ separates
$p(K)$ into two loops (that may intersect each other). Since $p(K)$ is
orientation preserving, either both of these loops are orientation preserving
or both are orientation reversing. 

If the two loops of $p(K)$ separated by
$d$ are orientation preserving, then the two points of $K$ over $d$ induce
the same orientation of the fiber. Thus under $\gamma_3$ the two branches of
$K$ over $d$  slide in the same direction and such double points $d$ do not
correspond to any input into $\Delta_{\wt\aslk}$.

If $d$ separates $p(K)$ into two orientation reversing loops, then the two
points of $K$ over $d$ induce the opposite orientations of the fiber over
$d$. Thus the two branches of $K$ over $d$ slide in the opposite directions
under the deformation $\gamma_3$, and such $d$ correspond to singular knots
arising under $\gamma_3$. However both loops of the arising singular knots
project to orientation reversing loops on $F$ and hence to orientation
reversing loops on $G$. Thus they are not contractible in $M$, and we
get that $\Delta_{\wt\aslk}(\gamma_3)=0$.

If one considers a framing of $K$ that is nowhere tangent to the fibers of
$p$, then it becomes clear that $\delta(\gamma_3)=0$. Since
$\delta(\gamma_1)=\Delta_{\wt\aslk}(\gamma_1)=0$ and
$\delta(\gamma_3)=\Delta_{\wt\aslk}(\gamma_3)=2$, we get that
$\Delta_{\wt\aslk}(\alpha^i)=\delta(\alpha^i)$. Since $i$ was taken to be nonzero we get
that $\Delta_{\wt\aslk}(\alpha)=\delta(\alpha)$.

{\bf Consider the case of $p(K)=1\in\pi_1(F)$.} Since $K\neq 1\in M$
we have that $K=f^k$, for some $k\neq 0$. Using~\ref{torussubstitution}
we can assume that $K$ is the $(1,k)$ toric knot on the torus that is $p^{-1}p(K)$. 

Clearly $p_*t(\alpha^2)$ is an orientation preserving loop on $F$. Let
$\gamma_{\alpha^2}$ be the isotopy of $K$ such that  
$p(\gamma_{\alpha^2}(x))$ is a small circle at every moment of time $x$, $\gamma_{\alpha^2}(x)$ is the $(1,k)$ toric knot
on the torus that is $p^{-1}p(\gamma_{\alpha_2}(x))$, and $t(\gamma_{\alpha^2})=t(\alpha^2)$. (It is
easy to verify that such an isotopy really does exist.) Since $\gamma_{\alpha^2}$ is an isotopy, we have $\Delta_{\wt\aslk}(\gamma_{\alpha^2})=0$.

Thus $\alpha^2=\gamma_{\alpha^2} \gamma_2^s$, for some $s\in\Z$. If one considers the
framing of $K$ such that the projections of the framing vectors to $F$ are
orthogonal to the circle $p(K)$, it is easy to verify that
$\delta(\gamma_{\alpha^2})=0$. Since $\Delta_{\wt\aslk}(\gamma_2)=\delta(\gamma_2)=2$ and
$\Delta_{\wt\aslk}(\gamma_{\alpha^2})=0$, we have $\Delta_{\wt\aslk}(\alpha^2)=\delta(\alpha^2)$. Thus
$\Delta_{\wt\aslk}(\alpha)=\delta(\alpha)$. This finishes the proof of the Lemma for all the cases.\qed
\end{emf}

\begin{lem}\label{vertical}
Let $q:M\rightarrow F$ be a compact oriented (not necessarily closed) 
Seifert-fibered manifold with $\pi_2(M)=0$. Let 
$\bar\mu_{\alpha}:S^1\times S^1\rightarrow M$ be homotopic to 
a vertical torus 
$\lambda:S^1\times S^1\rightarrow M$. Then $\Delta_{\wt\aslk}(\alpha)=\delta(\alpha)$.
\end{lem}

\begin{emf}{\em Proof of Lemma~\ref{vertical}.\/}
Let $M'\subset M$ be a thin neighborhood of $\Im(\lambda)$ such that it is
locally-trivially $S^1$-fibered over a thin neighborhood of
$q(\lambda)\subset F$. Let $q':M'\rightarrow F'$ be the corresponding
locally-trivial $S^1$-fibration. Since $M'$ is thin, we can assume that
$F'$ has nonempty boundary. Now Lemma follows from  Remark~\ref{torussubstitution}, $h$-principle~\ref{h-principle}, and Lemma~\ref{locallytrivial}.\qed
\end{emf}

{\em Thus we have reduced the Proof of Theorem~\ref{irreducible} 
to the case where $\bar\mu_{\alpha}$ is an essential torus 
inside of an oriented compact (not necessarily closed) 
Seifert-fibered space $M$ with $\pi_2(M)=0$ and $\pi_1(M)$ torsion free. 
(In particular $\pi_1(M)$ has to be infinite.) 
The arguments 
above give the proof of the Theorem in the case where 
$\bar\mu_{\alpha}$ is homotopic to a vertical  torus.\/}

\begin{emf}\label{Seiferthaken}
As it is shown in~\cite{Johannson}, see Propositions 5.13 (and remark after it) and 7.1, for most
Seifert-fibered manifolds all essential 
tori are homotopic to vertical ones for some
choice of Seifert-fibration structure.

The only compact Seifert-fibered manifolds where this statement is not proved 
in~\cite{Johannson} are: 
\begin{description}
\item[1] Seifert-fibered spaces with the orbit space $\R P^2$ and at most
one exceptional fiber; 
\item[2] Seifert-fibered spaces with $l$ exceptional fibers and 
the orbit space being $S^2$ with $m$ holes, $l+m\leq 3$.
\end{description}

If $M$ is Seifert-fibered over $\R P^2$ with at most one exceptional fiber,
then the orientation cover $S^2\rightarrow \R P^2$ induces a two fold cover 
$\tilde M\rightarrow M$. The manifold $\tilde M$ admits a structure of a
Seifert fibration over $S^2$ with at most two exceptional fibers. As it is
shown by Orlik~\cite{Orlik} such $\tilde M$ is a lens space. Since
$\pi_2(M)=0$, we get that $\pi_2(\tilde M)=0$. Thus $\pi_1(\tilde M)$ and
$\pi_1(M)$, are finite and we get the proof for this case.

Below we consider the case where $M$ is a Seifert-fibration over $S^2$ and
$l+m\leq 3$.  

{\bf Consider the case of $m\neq 0$.\/}
If $l=0$, then the proof follows from Lemma~\ref{locallytrivial}.

We shall use the standard presentation of $\pi_1$ of a Seifert fibered space, see for example~\cite{Jaco}. For a Seifert-fibered manifold $M$ over a sphere with $l\neq 0$ exceptional fibers and $m\neq 0$ holes we get that 
\begin{equation}\label{pi1seifert}
\pi_1(M)=\{c_1, \dots, c_l, d_1, \dots, d_{m-1}, f\big| c_jf=fc_j, d_jf=fd_j, c_j^{\alpha_j}=f^{\beta_j}\}; 
\end{equation}
where  $\alpha_j$ is the index of the $j$th exceptional fiber, $0<\beta_j<\alpha_j$, $f$ is the class of the regular fiber, and projections of $c_j$ and $d_j$ go around the $j$th cone point and $j$th hole on $S^2$, respectively.

For all the three cases $(l=1, m=1)$, $(l=1, m=2)$, and $(l=2, m=1)$ the image of the quotient homomorphism
$p:\pi_1(M)\to \pi_1(M)/ F$ by the normal subgroup $F$ generated by $f$ is the free product $\star$ of cyclic groups. 
%
%
%
%
%

An exercise in group theory shows that for any finitely generated groups $G_1, G_2$ 
and $g\in G_1\star G_2$: 
\begin{description}
\item[a] the centralizer of $g$ is isomorphic to the infinite 
cyclic group, provided that $g$ is not conjugate to an element of 
$G_1\star \{1\}<G_1\star G_2$ or to an element of $\{1\}\star G_2< G_1\star
G_2$;
\item[b] if $g\neq 1$ and there exists $\tilde g$ such that $\tilde g^{-1}g\tilde
g=g'\star\{1\}\in G_1\star G_2$, then the centralizer of $g$ is $\tilde g^{-1}
(Z(g')\star\{1\})\tilde g$, where $Z(g')$ is the centralizer of $g'\in G_1$.
\end{description}

Thus if $p(K)\neq 1$, then for every $\alpha\in \pi_1(\mathcal C)$ there exists $i\neq 0$ and $j\in\Z$ such that $p(t(\alpha^i))=p(t(\gamma_1^j))$, see~\ref{prooflocallytrivial}. If $p(K)=1$, then for every $\alpha\in \pi_1(\mathcal C)$ we get that $p(t(\alpha^2))=p(t(\gamma_{\alpha^2}))$. Using these observations and the fact that $f$ is in the center of $\pi_1(M)$ 
we get the proof in the cases $(l=1, m=1)$, $(l=1, m=2)$, and $(l=2, m=1)$ 
repeating the arguments of~\ref{prooflocallytrivial}.

{\bf Below we assume that $m=0$.\/}
If $m=l=0$, then since $M$ is irreducible
we have that $\pi_1(M)=\Z_i$, where $i\in\Z=H^2(S^2)$ is the Euler class of
the locally-trivial $S^1$-bundle $M\rightarrow S^2$. 
(If $i=0$, then $M=S^1\times S^2$.) 
Since in these cases $\pi_1(M)$ is finite, we get the proof.

If $m=0$ and $l=1$, or $m=0$ and $l=2$ then $M$ is a lens space, see 
for example Orlik~\cite{Orlik} p.~99. Since $M$ is irreducible it has a
finite fundamental group and we have finished the proof in these cases.

{\em Consider the case of $m=0$ and $l=3$.\/} 
Let $r,s,t\in\N$ be the multiplicities of the exceptional fibers.
The quotient group of $\pi_1(M)$ by the subgroup generated by
the regular fiber of $q$ is the triangle group $\Delta(r,s,t)$.

The result that can be found in~\cite{Orlik} pp.100--101 says that if
$\frac{1}{r}+\frac{1}{s}+\frac{1}{t}>1$, then $\pi_1(M)$ is finite and thus $M$
does not contain essential tori.  

The results of Hass, see~\cite{Hass} Theorem 1 and Lemma 2, (some of these
results
were independently obtained by Gao~\cite{Gao}) 
imply that every
essential mapping of a torus to a closed irreducible 
Seifert fibered manifold over an
orientable surface is homotopic to either a vertical immersed 
one or to a horizontal immersed
one. (A mapping is {\em horizontal\/} if 
it is everywhere transverse to the fibers of the fibration.) It is easy
to see from the work of Scott~\cite{Scott} 
that in the hyperbolic triangle group case, 
$\frac{1}{r}+\frac{1}{s}+\frac{1}{t}<1$, if the mapping of a torus is
horizontal then the torus has a negative curvature metric. Since this is
impossible, we get that every essential mapping of a torus in the hyperbolic
triangle group case is homotopic to a vertical one, and we get the proof for this case.

The only Euclidean triangle groups are $\Delta(2,3,6)$, $\Delta(2,4,4)$, and
$\Delta(3,3,3)$. Since we have already proved the Theorem for all the cases
where essential tori are homotopic to vertical tori, we get from the work of
Hass~\cite{Hass} that the only manifolds corresponding to Euclidean triangle
groups we have to consider are those that contain horizontal immersed tori.
Thus the Euler number of the Seifert-fibration should be zero, see for
example~\cite{Hatcher}.

Following the work of Kirk and Livingston~\cite{KirkLivingston} we observe  
that the only Seifert-fibered 
manifolds that correspond to the Euclidean triangle groups and have zero Euler
number are:
\begin{description}
\item[1] $M_{(2,3,6)}$ with Seifert invariants of the fibers
$\{(2,1), (3,-1), (6,-1)\}$;
\item[2] $M_{(2,4,4)}$ with Seifert invariants of the fibers $\{(2,1),
(4,-1), (4,-1)\}$; and
\item[3] $M_{(3,3,3)}$ with
Seifert invariants of the fibers $\{ (3,1), (3,1), (3,-2)\}$.
\end{description}
These spaces have structure of torus-bundles over a circle with a finite order monodromy
(see~\cite{KirkLivingston}, or~\cite{Hatcher} for a more detailed
explanation). 
The monodromy maps for $M_{(2,3,6)}$, $M_{(2,4,4)}$, and $M_{(3,3,3)}$ are
given respectively by the linear maps $\R^2/\Z^2\rightarrow \R^2/\Z^2$
determined by the matrices

\[
\begin{array}{lcr}
A=\left ( \begin{array}{lr}
0 & -1\\
1 &  1\\
\end{array}\right)
,&
B=\left ( \begin{array}{rr}
0 & 1\\
-1 & 0\\
\end{array}\right)
,&
C=\left( \begin{array}{lr}
0 & -1\\
1 & -1\\
\end{array}
\right ).
\\
\end{array}
\]

Thus their finite covering is $S^1\times S^1\times S^1$. To prove the Theorem for these three manifolds we will need the following Lemma.

\begin{lem}\label{coveringframingisotopy}
Let $q:S^1\times S^1\times S^1\rightarrow M$ be a finite covering.
Let $\wt {\mathcal C}$ be a connected component of the space of curves in $S^1\times S^1\times S^1$ and let $\tilde K \in \wt {\mathcal C}$ be a knot such that $K=q(\wt K)$ is a nonsingular knot. 
Put $\mathcal C=q (\wt {\mathcal C})$.
Then $\Delta_{\wt\aslk}(\alpha)=\delta(\alpha)$, for every $\alpha\in \pi_1(\mathcal C, K)$.
\end{lem}

\begin{emf}{\em Proof of Lemma~\ref{coveringframingisotopy}.\/}
The lifting of a loop $\alpha\subset \mathcal C$ 
to $\wt {\mathcal C}$ connects $\wt K$ to another lifting of $K$. 
Since the covering $q:S^1\times S^1\times S^1\rightarrow M$ is finite, there exists 
$i\neq 0$ such that $\alpha^i$ lifts to a homotopy $\wt {\alpha^i}$ of 
$\wt K$ to itself.

Similarly to~\ref{prooflocallytrivial} we get that $\wt {\alpha^i}=\gamma_{3,1}^{i_1}\gamma_{3,2}^{i_2}\gamma_{3,3}^{i_3}\gamma_2^m$, for some 
$i_1, i_2, i_3, m\in \Z$, where $\gamma_{3,1}, \gamma_{3,2}, \gamma_{3,3}$ are loops 
of type $\gamma_3$ with respect to the three obvious locally-trivial $S^1$-fibration structures on $S^1\times S^1\times S^1$. Put $\wt \gamma=\gamma_{3,1}^{i_1}\gamma_{3,2}^{i_2}\gamma_{3,3}^{i_3}\in \pi_1(\wt {\mathcal C})$ and $\gamma=q(\wt \gamma)\in \pi_1(\mathcal C)$. Clearly $\alpha^i=\gamma\gamma_2^m$.

Since $q$ is a covering and $\wt \gamma$ is an isotopy, $\gamma$ (which is not necessarily an isotopy) does not contain singular knots with a double point splitting the singular knot into two loops one of which is contractible. Thus $\Delta_{\wt \aslk}(\gamma)=0$. Clearly $\delta(\gamma)=0$, and by Lemma~\ref{betaforobstructions} $\delta(\gamma_2)=\Delta_{\wt \aslk}(\gamma_2)=2$. Thus $\delta(\alpha^i)=\Delta_{\wt \aslk}(\alpha^i)$, and since $i\neq 0$ we have $\delta(\alpha)=\Delta_{\wt \aslk}(\alpha)$. \qed
\end{emf}

For all the three manifolds $M\in\{M_{(2,3,6)}, M_{(2,4,4)},M_{(3,3,3)}\}$ 
the fundamental group of $M$
is a semi-direct product $\pi_1(T^2)\propto \pi_1(S^1)=
(\Z\oplus\Z)\propto\Z$. 
Let $m,l\in\pi_1(S^1\times S^1)$ be the classes of
respectively the meridian and the longitude, and let $f\in\pi_1(S^1)$ be the
generator. For $a,a'\in\Z\oplus\Z$ and $f^i, f^j\in\Z$ the product 
$(a, f^i)(a', f^j)\in (\Z\oplus\Z)\propto\Z$ is given by $(a \xi(f^i)(a'),
f^{i+j})$ where $\xi(f^i)(a')=f^ia'f^{-i}\in\pi_1(M)$. 

The action $\xi(f^k)(m^il^j)$ is calculated as follows.
Let $D\in M_2(\Z)$
be the monodromy matrix of the torus bundle $M\rightarrow S^1$. 
Put 
\[
\label{eqaction}
\left(\begin{array}{l}i' \\ j'\\
\end{array}
\right )= D^k\left(\begin{array}{l} i \\ j\\  
\end{array}
\right).
\]

Then 
\begin{equation}\label{equationaction}
\xi(f^k)(m^il^j)=m^{i'}l^{j'}.
\end{equation}

We prove the Theorem for $M=M_{(3,3,3)}$. The proof of the Theorem for 
$M\in\{M_{(2,3,6)}, M_{(2,4,4)}\}$ is completely analogous to the case of
$M=M_{(3,3,3)}$ but the calculations are a bit harder.

The matrix $C$ that describes the monodromy of $M_{(3,3,3)}\rightarrow S^1$
has order $3$. Thus $\xi(f^{3k})(a)=a$, for every $k\in\Z$ and
$a\in\pi_1(T^2)$.

Take a knot $K\in\mathcal C$ and $\alpha\in\pi_1(\mathcal C, K)$. Then 
$t(\alpha^3)=(a', f^{3k'})$ and $K=(a, f^{3k+n})$, for some
$a,a'\in\pi_1(T^2)$, $k,k'\in\Z$, and $n\in\{0,1,2\}$. Since $t(\alpha^3)$
commutes with $K$ in $\pi_1(M)$ we have that
\begin{equation}
(a\xi(f^{3k+n})(a'), f^{3k+3k'+n})=(a'\xi(f^{3k'})(a), f^{3k+3k'+n}).
\end{equation} 
Since $\pi_1(T^2)$ is commutative and $\xi(f^{3k'})(a)=a$, we get that 
\begin{equation}\label{fixedelement}
\xi(f^{3k+n})(a')=a'.
\end{equation}

{\em Consider the case of $n\in\{1,2\}$.\/}
Then one uses~\eqref{equationaction} to verify that 
$a'=1=0\oplus0\in\Z\oplus\Z$ (provided that $n\in\{1,2\}$). Thus
$t(\alpha^3)=(1,f^{3k'})\in(\Z\oplus\Z)\propto \Z$.
A straightforward calculation shows that $K^3=(a, f^{3k+n})^3=(1,
f^{3(3k+n)})\in\pi_1(M_{(3,3,3)}, K(1))$, thus $K^{3k'}=(1, f^{3k'(3k+n)})=
(1, f^{3k'})^{3k+n}=t(\alpha^{3(3k+n)})\in\pi_1(M_{(3,3,3)}, K(1))$.

Let $\gamma_1\in\pi_1(\mathcal C, K)$ be the loop introduced
in~\ref{gamma2}.
Then $t(\gamma_1)=K$ and
$t(\gamma_1^{3k'})=t(\alpha^{3(3k+n)})$. Thus
$\alpha^{3(3k+n)}=\gamma_1^{3k'}\gamma_2^s$, for some $s\in\Z$. 
Since $\Delta_{\wt\aslk}(\gamma_1)=\delta(\gamma_1)=0$ and
$\Delta_{\wt\aslk}(\gamma_2)=\delta(\gamma_2)=2$, we get that
$\Delta_{\wt\aslk}(\alpha^{3(3k+n)})=\delta(\alpha^{3(3k+n)})$. Since $3(3k+n)\neq 0$,
we have $\Delta_{\wt\aslk}(\alpha)=\delta(\alpha)$ and this finishes the proof for
$M=M_{(3,3,3)}$ and $K=(a, f^{3k+n})\in\pi_1(M)$ with $n\in\{1,2\}$.

{\em Consider the case where $n=0$.\/} 
Then $K$ is liftable to the total space of
the three fold covering $S^1\times S^1\times S^1\rightarrow M_{(3,3,3)}$, and  Lemma~\ref{coveringframingisotopy} implies the statement of the Theorem.

This finishes the proof of Theorem~\ref{irreducible} 
for $M=M_{(3,3,3)}$, and thus for all the
remaining cases. This is also the end of the Proof of Theorem~\ref{aslk}.\qed
\end{emf}

\begin{rem}\label{essentialorientable}
Most parts of the proof  above work if one tries to construct the
$\aslk$ invariant for knots homotopic into an irreducible summand $M_j$ of a closed oriented manifold $M\neq (S^1\times S^2)\#M'$. 
There are exactly two steps that collapse. They are: 
\begin{description} 
\item[1] The step in the proof of Lemma~\ref{locallytrivial} when $p(K)$ has double points $d$ that split $p(K)$ into two orientation reversing loops. In this case it is easy to construct examples showing that $\aslk$ does not always exist.
However, if for $K$ such $d$ do not exist (or if the total input of all of them into $\Delta_{\aslk}(\gamma_3)$ is zero) this step works for $\aslk$.
\item[2] The proof of Lemma~\ref{coveringframingisotopy} does not work for the $\aslk$ invariant. Thus we can not prove the existence of the $\aslk$ for  
$M_j\in\{M_{(2,3,6)}, M_{(2,4,4)},M_{(3,3,3)}\}$. 
\end{description}

In particular the $\aslk$ invariant does exist 
if the irreducible summand $M_j$ of $M$ (into which $K$ is homotopic) is not one 
of $M_{(2,3,6)}, M_{(2,4,4)},M_{(3,3,3)}$ and 
the components of the characteristic submanifold of $M_j$ do not admit a structure of a Seifert-fibration over a nonorientable surface. 
\end{rem}

\section{Proof of Theorem~\ref{framed}}\label{proofframed}
Let $K_f$ be $K$ with some framing. Put $K_f^i$, $i\in \Z$, to be framed knots defined in~\ref{Introduction}.
If $|K|\neq \infty$, then there exists $i\neq j\in\Z$ such that $K_f^i$ and $K_f^j$ are isotopic framed knots.

Since the isotopy that changes $K_f^i$ to $K_f^j$ happens in a compact part of $M$, and $M$ is realizable as a connected sum $(S^1\times S^2)\#M$ if and only if $M$ contains a non-separating sphere, we can assume in the proof that $M$ is compact.

The following result~\ref{Matveevlemma} was told to the author by
S.~Matveev~\cite{Matveev} and it allows us to assume that $M$ is closed.

\begin{lem}\label{Matveevlemma}
Let $M$ be a connected compact oriented $3$-manifold that does not contain
an embedded non-separating $2$-sphere, then $M$ is realizable as a submanifold
of a closed oriented $3$-manifold $N$ such that $N$ also does not contain an
embedded non-separating $2$-sphere.
\end{lem}

\begin{emf}{\em Proof of Lemma~\ref{Matveevlemma}.\/}
Let $D^2_i$, $i\in I$, be a maximal collection of 
embedded non-intersecting 
disks in $M$ with the properties that 
$\p D^2_i \subset \p M$ and the union of $D^2_i$ does not
separate $M$. Let $M'$ be a connected compact 
manifold obtained by cutting $M$ along all of the disks $D^2_i$, $i\in I$. 
Clearly there are no embedded non-separating disks $D^2$ in $M'$ with $\p
D^2\subset \p M'$. Let $N$ be a closed oriented $3$-manifold that is the double
of $M'$ i.e. $N$ is obtained by gluing together two copies $M'_1$ and
$M'_2$ of $M'$ with
the opposite orientations along the identity automorphism of the boundary.
Let $\tilde M$ be an oriented submanifold of $N$ obtained 
by gluing for every $i\in I$ 
a thin solid tube $D^2\times [0,1]$ to $M'_1$ 
inside of $M'_2$ so that $D^2\times \{0\}$ and $D^2\times \{1\}$ are glued to
the two copies of $D^2_i\subset
\p M'_1$; $D^2\times(0,1)$ is mapped into $N\setminus M'_1$; 
and the thin tubes that correspond to different $i_1, i_2\in I$ 
do not intersect pairwise.
Clearly $\tilde M$ is diffeomorphic to $M$. 

Below we show that $N$ does not contain embedded non-separating $2$-spheres. Let
$D^2_j$, $j\in J$, be a maximal collection of embedded pairwise
non-intersecting $2$-disks in $M'$ with the properties 
that $\p D^2_j\subset \p M'$ and none of the pieces
they cut $M'$ into is $D^2\times [0,1]$. Let $S^2_j$, $j\in J$, be the embedded 
$2$-spheres in $N$ obtained by gluing the two copies of $D^2_j$ located in the
two copies of $M'$ along the common boundary. The spheres $S^2_j$ give a
decomposition of  $N$ into a connected sum of closed $3$-manifolds 
$N_k$, $k\in K$, such that every $N_k$ is a double of a compact
$3$-manifold $M_k$ with the property that there are no embedded 
$2$-disks $D^2$ in $M_k$
with $\p D^2\subset \p M_k$ and $\p D^2$ being a non-contractible loop in $\p
M_k$. The manifolds $M_k$ do not contain non-separating embedded spheres by the
assumption of the Lemma, and the statement of the Lemma follows from the
following Proposition.
\qed
\end{emf}
 
\begin{prop}\label{double}
Let $M$ be a compact connected oriented $3$-manifold that does not contain embedded
non-separating $2$-spheres and embedded $2$-disks $D^2$ with $\p
D^2\subset \p M$ and $\p D^2\neq 1\in \pi_1(\p M)$. 
Then the double $N$ of $M$ does not contain embedded
non-separating $2$-spheres.
\end{prop}

\begin{emf}{\em Proof of Proposition~\ref{double}.\/}
Assume that $N$ contains an embedded non-separating $2$-sphere. 
Take such a sphere $i:S^2\rightarrow N$ that is transverse to $\p M$ and 
that has the minimal number $c$ of connected components of
the intersection $i(S^2)\cap \p M$. We show that $c=0$ and get a
contradiction with the assumptions of the Proposition. 

Assume that $c>0$. Take a circle $\gamma$ of $i(S^2)\cap \p M$ that bounds an
embedded disk in $M$. By assumptions of the proposition $\gamma$ bounds a
disk $D^2\subset \p M$. Take one of the inner most circles $\tilde \gamma$ of 
$D^2\cap i(S^2)$. 
Then $\tilde \gamma$ bounds an embedded disk $\tilde D^2\subset \p M$ such
that $\Int \tilde D^2\cap i(S^2)=\emptyset$. 
Cut $i(S^2)$ along $\tilde \gamma$, take the two connected
components of $i(S^2)\setminus \tilde \gamma$ and glue the boundaries of
them with the copies of $\tilde D^2$. We get two embedded spheres.
Push slightly the glued disks off $\p (M)$ into the corresponding copies of
$M$ and get that for both embedded
spheres the number of connected components of their intersection with $\p M$
is less than $c$. One of these embedded 
spheres is non-separating for  homological reasons. This implies that $c=0$.\qed
\end{emf}

\begin{emf}\label{framedcomponents}
Now we have to prove Theorem~\ref{framed} in the case where $M$ is closed.

Every oriented $3$-dimensional manifold $M$ is parallelizable, and hence it
admits a $\spin$-structure. A framed  
curve $ K_f$ in $M$ represents a loop in the principal $SO(3)$-bundle of $TM$.
The $3$-frame corresponding to a point of $ K_f$ is the velocity
vector, the framing vector, and the unique third vector of unit length such
that the $3$-frame defines the positive orientation of $M$.
Clearly the values of a spin structure on the loops in the principal
$SO(3)$-bundle of $TM$ that correspond to 
$K_f^i$ and $K_f^j$ are different, provided that $i-j$ is odd.
Since the value of a spin structure depends only on the connected component
of the space of framed curves, we get that  $K_f^i$ and
$K_f^j$ with $i-j$ odd are not homotopic. Hence in our case $i-j=2k$ is even.

For simplicity we assume that $k>0$. Consider a homotopy of $K_f^j$ that is 
a framed version of $\gamma_2^k$, see~\ref{gamma2} for the definition of $\gamma_2$. In the end of this homotopy we get $K_f^i$. 

Now we use Theorem~\ref{aslk}. In the case where $K$ is homotopic to a loop in an irreducible summand of $M=\#_{j\in J}M_j$ by Theorem~\ref{aslk} the $\wt \aslk$ invariant is well-defined. Since one of the loops of the singular knot arising under $\gamma_2$ is contractible, we get that $\wt \aslk(K_f^j)-\wt \aslk(K_f^i)=2k$. Thus every homotopy between $K_f^i$ and $K_f^j$ involves at least $k$ passages through a double point and 
$K_f^i$ is not isotopic to $K_f^j$.

If $K$ is not homotopic to a loop in one 
of the irreducible summands of $M=\#_{j\in J}M_j$, then the proof is obtained in a similar way using the $\aslk$ invariant that is known to exist by Theorem~\ref{aslk}. \qed
\end{emf}

{\bf Acknowledgments.}
I am very grateful to I.~Fedorova, Yu.~Rudyak, and A.~Vylejneva for many 
enlightening discussions.  
I am very grateful to S.~Matveev for the discussions about
non-separating spheres in $3$-manifolds and telling me the proof of Lemma~\ref{Matveevlemma}.
I am thankful to N.~A'Campo, A.~Cattaneo,
P.~Kirk, C.~Livingston, and P.~Scott
for the useful mathematical discussions and suggestions.

The ideas of some parts of the proof are inspired by the works of
V.~Arnold~\cite{Arnoldcurves}, E.~Kalfagianni~\cite{Kalfagianni}, and 
P.~Kirk and C.~Livingston~\cite{KirkLivingston}
and~\cite{KirkLivingstontype1invariants}. 

This paper was written during my stay at the Institute for Mathematics,
Z\"urich University. It was substantially revised at Dartmouth College, when I was supported by the Dartmouth free term research salary.
I  thank Z\"urich University and Dartmouth College for the excellent working conditions.

\end{document}